\documentclass{article}
\usepackage{amssymb}

\newenvironment{proof}[1][Proof]{\noindent\textbf{#1.} }{\ \rule{0.5em}{0.5em}}
\input{tcilatex}

\begin{document}

\begin{center}
\begin{eqnarray*}
&&\text{{\LARGE \ }{\large On the distances between Pisot numbers generating}%
} \\
&&\text{{\large \ \ \ \ \ \ \ \ \ \ \ \ \ \ \ \ \ \ \ \  the same number
field }}
\end{eqnarray*}

\smallskip

by\ Toufik Za\"{\i}mi
\end{center}

\textbf{Abstract. \ }\textit{A well-known result, due to Meyer, states that
the set }$\wp _{K}$\textit{\ of Pisot numbers, generating a real algebraic
number field }$K,$\textit{\ is uniformly discrete and relatively dense in
the interval }$[1,\infty ).$\textit{\ In the present paper, we show that }$%
\wp _{K}\cup \{1\}\subset \wp _{K}-\wp _{K}$\textit{\ \ and the complement
of }$\wp _{K}$ \textit{in }$\wp _{K}-\wp _{K}$\textit{\ is not finite. Also,
we prove\ that if }$K$\ \textit{is totally real, then the elements of }$\wp
_{K}-\wp _{K}$ \textit{are the algebraic integers of }$K$\textit{\ whose
images under the action of all embeddings of }$K$\textit{\ into }$\mathbb{R}%
, $\textit{\ other than the identity of }$K,$\textit{\ belong to the
interval }$(-2,2),$ \textit{and so any Salem trace number }$\beta $\textit{\
may be written as a difference of two Pisot numbers generating the field} $%
\mathbb{Q}(\beta ).$

\medskip

\textbf{1. Introduction}

\textbf{\ }A Pisot number is a real algebraic integer greater than $1$ whose
other conjugates are of modulus less than $1,$ and the set of such numbers
is traditionally denoted by $S.$ Let throughout $K$ denote a real algebraic
number field of degree $d,$ and let 
\[
\wp _{K}:=\{\theta \mid \theta \in S\cap K,\text{ \ }K=\mathbb{Q}(\theta
)\}. 
\]%
Since there are at most a finite number of algebraic integers of fixed
degree and having all conjugates in a bounded subset of the complex plane,
the set of Pisot numbers of a given degree is a discrete subset of $%
[1,\infty ),$ i. e., \ any finite subinterval of $[1,\infty )$ contains at
most a finite number of such numbers, and so the family $\wp _{K}$ is also
discrete in $[1,\infty ).$

Pisot [13] was the first to investigate the set $\wp _{K}$ and he showed, in
particular, that $\wp _{K}$ contains units, whenever $d\geq 2.$ Fan and
Schmeling [7] have proved that $\wp _{K}$ is relatively dense in $[1,\infty
),$ that is, there is a real number $\rho >0$ such that every subinterval of 
$[1,\infty )$ of the form $(\varepsilon ,\varepsilon +\rho ]$ contains at
least one element of $\wp _{K};$ to be more precise, we say in this case
that $\wp _{K}$ is $\rho $-dense in $[1,\infty ).$

Using some theorems, due to Meyer, on harmonious sets (see for instance [8],
[10], and [12]), the author pointed out, in [19], that the family $\wp
_{K}^{\prime }:=\wp _{K}\cup (-\wp _{K})$ is a real Meyer set, i. e., there
is a finite subset $F$ of the real line $\mathbb{R}$ such that $\wp
_{K}^{\prime }-\wp _{K}^{\prime }\subset \wp _{K}^{\prime }+F,$ and a
generalization of this result for complex Pisot numbers was given in [2]. A
real Meyer set is, in particular, a relatively dense subset of $\mathbb{R}$,
which is also uniformly discrete, that is, there is a positive real number $%
\rho $ such that each real interval of the form $(\varepsilon ,\varepsilon
+\rho ]$ contains at most one element of this set [12].

It is worth noting that the set of Pisot numbers of a given degree $d\geq 2,$%
\ is not uniformly discrete. Indeed, it easy to see, by Rouche's theorem,
that a polynomial of the form $x^{d}-nx^{d-1}-k,$ where $k\in \{1,2\}$ and $%
n\in \lbrack 4,\infty )\cap \mathbb{N},$ is the minimal polynomial of a
Pisot number lying in $(n,n+k/n),$ and so the interval $(n,n+2/n)$ contains
two distinct Pisot numbers of degree $d.$

Let $\theta _{1}<\theta _{2}<\theta _{3}<\cdot \cdot \cdot $ \ denote the
elements of $\wp _{K},$ 
\[
D_{K}:=\{\theta _{m}-\theta _{n}\mid (m,n)\in \mathbb{N}^{2},\text{ }n>m\}, 
\]%
and

\[
\tciFourier _{K}:=\{\theta _{n+1}-\theta _{n}\mid \text{ }n\in \mathbb{N}\}. 
\]%
Then, $\wp _{K}-\wp _{K}=(-D_{K})\cup \{0\}\cup D_{K}$ and 
\begin{equation}
\tciFourier _{K}\subset D_{K}\subset E_{K}:=\{\beta \in \mathbb{Z}_{K}\cap
(0,\infty )\mid (\left\vert \sigma _{2}(\beta )\right\vert ,...,\left\vert
\sigma _{d}(\beta )\right\vert )\in (0,2)^{d-1}\},  \tag{1}
\end{equation}%
\bigskip where $\mathbb{Z}_{K}$ is the ring of the integers of $K,$ and $%
\sigma _{2},...,\sigma _{d}$ are the distinct embeddings of $K\mathbb{\ }$%
into $\mathbb{C}$ other than the identity of $K.$

Because the set $\wp _{K}$ is $\rho $-dense in $[1,\infty ),$ for some $\rho
>0,$ we have $(\theta _{n},\theta _{n}+\rho ]\cap \wp _{K}\neq \varnothing $
\ for each $n,$ and so $\theta _{n+1}-\theta _{n}\leq \rho .$ Further, as $%
\tciFourier _{K}\subset E_{K},$ the degree and the absolute values of the
conjugates of each element of $\tciFourier _{K}$ are, respectively, majored
by $d$ and $\max (\rho ,2),$ and hence the set $\tciFourier _{K}$ is finite.

Let%
\[
d_{1}<\cdot \cdot \cdot <d_{\func{card}(\tciFourier _{K})} 
\]%
be the elements of $\tciFourier _{K}.$ Then, $\wp _{K}$ is $\max \{\theta
_{1}-1,$ $d_{\func{card}(\tciFourier _{K})}\}$-dense in $[1,\infty ),$ and $%
\max \{\theta _{1}-1,$ $d_{\func{card}(\tciFourier _{K})}\}$ is the infimum
of those numbers $\rho $ for which $\wp _{K}$ is $\rho $-dense in $[1,\infty
).$ Also, the distance between any two distinct elements of $\wp _{K}$ is at
least $d_{1}$ (this implies that $\wp _{K}$ is uniformly discrete), and each
element of $D_{K}$ may be written 
\begin{equation}
\theta _{n}-\theta _{m}=(\theta _{n}-\theta _{n-1})+(\theta _{n-1}-\theta
_{n-2})+\cdot \cdot \cdot +(\theta _{m+1}-\theta _{m})=l_{1}d_{1}+\cdot
\cdot \cdot +l_{\func{card}(\tciFourier _{K})}d_{\func{card}(\tciFourier
_{K})},  \tag{2}
\end{equation}%
for some non-negative integers $l_{1},$ $...,$ $l_{\func{card}(\tciFourier
_{K})},$ not all zero.

For example, if $K=\mathbb{Q},$ then $\wp _{\mathbb{Q}}=\mathbb{N\diagdown \{%
}1\},$ $\wp _{\mathbb{Q}}-\wp _{\mathbb{Q}}=\mathbb{Z},$ $D_{\mathbb{Q}}=%
\mathbb{N}=E_{\mathbb{Q}}$ and $\tciFourier _{\mathbb{Q}}=\{1\}.$ It is of
interest to determine the elements of the sets $\tciFourier _{K}$ and $%
D_{K}. $

\medskip

\textbf{Theorem 1.1 } \textit{With the notation above, we have the following.%
}

\textit{(i) }$1\in D_{K};$\textit{\ \ }

\textit{(ii) }$\wp _{K}\varsubsetneq D_{K};$

\textit{(iii) if }$K\neq \mathbb{Q},$\textit{\ then the complement of }$\wp
_{K}$\textit{\ in }$D_{K}$\textit{\ is not finite, and }$\func{card}%
(\tciFourier _{K})\geq 2.$

\bigskip

Theorem 1.1(i) implies the result below, shown recently by Dubickas in [5].

\bigskip

\textbf{Corollary 1.2 [5] } \textit{If }$K=\mathbb{Q}(\tau )$\textit{\ for
some Salem number }$\tau ,$\textit{\ then }$\tau \in D_{K}.$

\medskip

Recall that a Salem number is a real algebraic integer $\tau >1$ whose other
conjugates are of modulus \ at most $1,$ with at least one conjugate of
modulus $1.$ Then, $\tau $ has two real conjugates, namely $\tau $ and $\tau
^{-1},$ and $\deg (\tau )-2\geq 2$ conjugates, say $\tau _{2}^{\pm
1},...,\tau _{d}^{\pm 1},$ lying on the unit circle. Also, the algebraic
integer $\tau +\tau ^{-1}>2,$ called a Salem trace number [20], is of degree 
$d,$ and its other conjugates are the numbers $\tau _{2}+\tau _{2}^{-1},$ $%
...,$ $\tau _{d}+\tau _{d}^{-1},$ lying in the interval $(-2,2).$

\ Conversely, an algebraic integer $\beta >2$ of degree $d\geq 2,$ whose
other conjugates belong the interval $(-2,2),$ is a Salem trace number
associated to some Salem number $\tau $ of degree $2d,$ via the relation $%
\beta =\tau +\tau ^{-1}$ (for more details, see the proofs of Lemmas 3.1 and
3.3).

By analogy with Salem trace numbers, we call a positive algebraic integer of
the form $u:=e^{i2k\pi /n}+e^{-i2k\pi /n}=2\cos (2k\pi /n),$ where $n\in 
\mathbb{N}\cap \lbrack 4,\infty ),$ $k\in \{1,...,n-1\}$ and $\gcd (k,n)=1,$
a root of unity trace number. Then, $\deg (u)=\varphi (n)/2,$ where $\varphi 
$ is the Euler totient function, and the conjugates of $u$ are the numbers $%
2\cos (2l\pi /n)\in (-2,2),$ where $l\in $ $\{1,...,[n/2]\},$ $[.]$ is the
integer part function, and $\gcd (l,n)=1.$

Clearly, the set, say $U_{K},$ of root of unity trace numbers lying in $K$
is finite, and $1=2\cos (2\pi /6)\in U_{K};$ we denote by 
\[
u_{1}<\cdot \cdot \cdot <u_{\func{card}(U_{K})} 
\]%
the elements of $U_{K}.$

For the case where the field $K$ is totally real, i.e., the conjugates of $K$
are all real, the set of Salem trace numbers generating $K$ over $\mathbb{Q}%
, $ and the set $U_{K}$ form a partition of $D_{K},$ as stated by the
following result.

\bigskip

\textbf{Theorem 1.3} \textit{Let }$K$\textit{\ be a totally real number
field of degree }$d\geq 2.$\textit{\ Then, the assertions below are true.}

\textit{(i) }$D_{K}=E_{K};$

\textit{(ii) }$U_{K}$\textit{\ }$=E_{K}\cap (0,2)=D_{K}\cap (0,2);$

\textit{(iii) the sets }$T_{K}:=E_{K}\cap (2,\infty )=D_{K}\cap (2,\infty )$%
\textit{\ and} $U_{K}$ \textit{form a partition of }$D_{K},$ \textit{and the}
\textit{degree of each element of }$T_{K}$\textit{\ is equal to} $d;$

\textit{(iv)} $\func{card}(\tciFourier _{K})\geq 2^{d-1};$

\textit{(v) }$\tciFourier _{K}\cap (0,2)\subset U_{K},$ $d_{1}=u_{1}\leq 1,$ 
\textit{and} $d_{2}=u_{2}$ \textit{whenever} $\func{card}(U_{K})\geq 2;$

\textit{(vi) if} $\func{card}(U_{K})=1,$ \textit{then }$d_{1}=1$ \textit{and}
$d_{2}=\min T_{K}$\textit{\ is the smallest Salem trace number of degree }$%
d, $\textit{\ lying in} $K.$

\bigskip

Notice that Theorem 1.3(vi) holds, for example, when $d$ is prime and the
discriminant of the field $K$ is sufficiently large or when $d$ is prime and 
$2d$ is not a totient number, so that the field $K$ cannot be embedded in a
cyclotomic field.

In the next section, we determine explicitly the elements of $\wp _{K},$ $%
\tciFourier _{K}$ and $U_{K},$ when $K$ is a real quadratic field. The
proofs of Theorem 1.1 and its corollary are given in Section 4, and the
proof of Theorem 1.3 is postponed to the last section. These proofs are
based on some lemmas presented in Section 3.

Throughout, when we speak about conjugates, the norm and the degree of an
algebraic number, without mentioning the basic field, this is meant over $%
\mathbb{Q}.$ Also, the degree, the discriminant and the conjugates of a
number field are considered over $\mathbb{Q},$ and an integer means a
rational integer. \ All computations are done using the system Pari [15].

\bigskip

\bigskip

\textbf{2. The quadratic case}

Let $K=\mathbb{Q}(\sqrt{m})$ be a real quadratic field, where $m\equiv 2,3$ $%
\func{mod}4$ is a square-free integer, and let $\theta \in \wp _{K}.$ Then, $%
\theta =a+b\sqrt{m}$ for some $(a,b)\in \mathbb{Z}\times \mathbb{Z\diagdown }%
\{0\},$ and the inequalities $a+b\sqrt{m}>1$ and $\left\vert a-b\sqrt{m}%
\right\vert <1$ yield $b\geq 1,$ $a\in \{\left\lfloor b\sqrt{m}\right\rfloor
,\left\lfloor b\sqrt{m}\right\rfloor +1\},$ and so 
\[
\theta =\alpha _{b}:=\left\lfloor b\sqrt{m}\right\rfloor +b\sqrt{m}\text{ \
\ or \ \ }\theta =\widehat{\alpha }_{b}:=1+\alpha _{b}\}. 
\]%
It is clear that for each $b\in \mathbb{N},$ $\{\alpha _{b},\widehat{\alpha }%
_{b}\}\subset \wp _{K},$ $\widehat{\alpha }_{b}-\alpha _{b}=1,$ and $\alpha
_{b+1}-\widehat{\alpha }_{b}=\left\lfloor (b+1)\sqrt{m}\right\rfloor +(b+1)%
\sqrt{m}-(1+\left\lfloor b\sqrt{m}\right\rfloor +b\sqrt{m})=\left\lfloor
(b+1)\sqrt{m}\right\rfloor +\sqrt{m}-(1+\left\lfloor b\sqrt{m}\right\rfloor
).$ Because $\left\lfloor b\sqrt{m}\right\rfloor +\left\lfloor \sqrt{m}%
\right\rfloor <(b+1)\sqrt{m}<\left\lfloor b\sqrt{m}\right\rfloor
+\left\lfloor \sqrt{m}\right\rfloor +2,$ we see that $\left\lfloor (b+1)%
\sqrt{m}\right\rfloor -\left\lfloor b\sqrt{m}\right\rfloor \in
\{\left\lfloor \sqrt{m}\right\rfloor ,\left\lfloor \sqrt{m}\right\rfloor
+1\} $ and 
\[
\alpha _{b+1}-\widehat{\alpha }_{b}\in \{\left\lfloor \sqrt{m}\right\rfloor
-1+\sqrt{m},\text{ }\left\lfloor \sqrt{m}\right\rfloor +\sqrt{m}\}. 
\]%
Therefore, the elements of $\wp _{K}$ may be labelled as follows:%
\[
\theta _{1}:=\alpha _{1}<\theta _{2}:=\widehat{\alpha }_{1}<\theta
_{3}:=\alpha _{2}<\theta _{4}:=\widehat{\alpha }_{2}<\theta _{5}:=\alpha
_{3}<\cdot \cdot \cdot 
\]%
Since the sequence $(b\sqrt{m})_{b\in \mathbb{N}}$ is dense modulo one (this
also follows from Lemma 3.1 with $n=1$), we get from the equivalences $%
\alpha _{b+1}-\widehat{\alpha }_{b}=\left\lfloor \sqrt{m}\right\rfloor -1+%
\sqrt{m}\Leftrightarrow b\sqrt{m}-\left\lfloor b\sqrt{m}\right\rfloor <1+$ $%
\left\lfloor \sqrt{m}\right\rfloor -\sqrt{m}$ and $\ \alpha _{b+1}-\widehat{%
\alpha }_{b}=\left\lfloor \sqrt{m}\right\rfloor +\sqrt{m}\Leftrightarrow b%
\sqrt{m}-\left\lfloor b\sqrt{m}\right\rfloor >1+\left\lfloor \sqrt{m}%
\right\rfloor -\sqrt{m}$ that 
\[
\{\alpha _{b+1}-\widehat{\alpha }_{b}\mid b\in \mathbb{N}\}=\{\left\lfloor 
\sqrt{m}\right\rfloor -1+\sqrt{m},\text{ }\left\lfloor \sqrt{m}\right\rfloor
+\sqrt{m}\}, 
\]%
and so 
\[
\tciFourier _{K}=\{1,\text{ }\left\lfloor \sqrt{m}\right\rfloor -1+\sqrt{m},%
\text{ }\left\lfloor \sqrt{m}\right\rfloor +\sqrt{m}\}. 
\]%
Finally, notice that if $a+b\sqrt{m}\in $ $U_{K}$ for some $(a,b)\in \mathbb{%
Z}^{2},$ then $a+b\sqrt{m}\in (0,2),$ $a-b\sqrt{m}\in (-2,2),$ and so $%
(a,b)=(1,0)$ (resp. $(a,b)\in \{(1,0),(0,1)\})$ when $m\geq 6$ (resp. when $%
m\in \{2,3\}).$ Therefore, we have the following.

\medskip \medskip

\textbf{Proposition 2.1} \textit{Let }$K=\mathbb{Q}(\sqrt{m})$\textit{\ be a
real quadratic field, where }$m\equiv 2,3$\textit{\ }$\func{mod}4$\textit{\
is a square-free integer. If} $m=2$ (\textit{resp.} $\ m=3,$ $m\geq 6),$ 
\textit{then}%
\[
U_{K}=\{1,\sqrt{2}=2\cos (2\pi /8)\}\subset \tciFourier _{K}=\{1,\sqrt{2},%
\text{ }1+\sqrt{2}=\min T_{K}\} 
\]%
\textit{\ \ }(\textit{resp. }%
\[
U_{K}=\{1,\sqrt{3}=2\cos (2\pi /12)\}\subset \tciFourier _{K}=\{1,\sqrt{3},1+%
\sqrt{3}=\min T_{K}\}, 
\]%
\textit{\ }%
\[
U_{K}=\{1\}\subset \tciFourier _{K}=\{1,\text{ }\left\lfloor \sqrt{m}%
\right\rfloor -1+\sqrt{m}=\min T_{K},\text{ }1+\min T_{K}\}). 
\]%
In fact, the equality $\min T_{K}=\left\lfloor \sqrt{m}\right\rfloor -1+%
\sqrt{m}$\textit{\ }(resp. \textit{\ }$\min T_{K}=1+\sqrt{2},$ $\min T_{K}=1+%
\sqrt{3}),$ for $m\geq 6$ \ (resp. for $\ m=2,$ for $m=3)$ follows from
Theorem 1.3(vi) (resp. from a simple computation).

In the same way, we obtain the following assertion.

\bigskip

\textbf{Proposition 2.2} \textit{Let }$K=\mathbb{Q}(\sqrt{m})$\textit{\ be a
real quadratic field, where }$m\equiv 1$\textit{\ }$\func{mod}4$\textit{\ is
a square-free integer. Then,} 
\[
U_{K}=\{1\}\subset \tciFourier _{K}=\{1,\text{ }\frac{-3+\left\lfloor \sqrt{m%
}\right\rfloor +\sqrt{m}}{2}=\min T_{K},\text{ }1+\min T_{K}\} 
\]%
\textit{when }$\left\lfloor \sqrt{m}\right\rfloor $\textit{\ is an even
integer greater than }$2,$ 
\[
U_{\mathbb{Q}(\sqrt{5})}=\tciFourier _{\mathbb{Q}(\sqrt{5})}=\{\frac{-1+%
\sqrt{5}}{2}=2\cos (\frac{2\pi }{5}),\text{ }1,\text{ }\frac{1+\sqrt{5}}{2}%
=2\cos (\frac{2\pi }{10})\}, 
\]%
\textit{and} 
\[
U_{K}=\{1\}\subset \tciFourier _{K}=\{1,\text{ }\frac{-2+\left\lfloor \sqrt{m%
}\right\rfloor +\sqrt{m}}{2}=\min T_{K},1+\min T_{K}\} 
\]%
\textit{when }$\left\lfloor \sqrt{m}\right\rfloor $\textit{\ is odd.}

\medskip

Finally, notice that a short calculation gives that $\min T_{\mathbb{Q}(%
\sqrt{5})}=(3+\sqrt{5})/2,$ and it follows, by Propositions 2.1 and 2.2,
that 
\[
\min_{K\text{ real quadratic field }}(\min T_{K})=\frac{1+\sqrt{13}}{2}; 
\]%
thus the smallest quartic Salem number $\tau _{0}$ satisfies the
(well-known) equation $\tau _{0}+1/\tau _{0}=(1+\sqrt{13})/2.$

\bigskip

\textbf{Remark 2.3} Let $K$ be a totally real cubic field and let $\alpha
\in U_{K}.$ Then, $\alpha =1$ or $\mathbb{Q}(\alpha )=K.$ Suppose $\alpha $
cubic. Then, $\alpha $ is a conjugate of $2\cos (\frac{2\pi }{n}),$ where $%
\varphi (n)=6.$ Hence, $n\in \{7,14,9,18\},$ and if $n\in \{7,14\}$ (resp. $%
n\in \{9,18\}),$ then the discriminant $\func{disc}(K)$ of $K$ is equal to $%
49$ and 
\[
U_{K}=\{-2\cos (\frac{4\pi }{7}),1,2\cos (\frac{2\pi }{7}),-2\cos (\frac{%
6\pi }{7})\} 
\]
(resp. \ is equal to $81$ and 
\[
U_{K}=\{2\cos (\frac{4\pi }{9}),1,2\cos (\frac{2\pi }{9}),-2\cos (\frac{8\pi 
}{9})\}). 
\]
Also, if $\func{disc}(K)\notin \{49,81\},$ then $U_{K}=\{1\}.$

Since the smallest Salem number of degree $6,$ say $\tau ,$ is given by the
equality $\tau ^{6}-\tau ^{4}-\tau ^{3}-\tau ^{2}+1=0,$ the smallest cubic
Salem trace number $\beta =\tau +1/\tau $ satisfies $\beta ^{3}-4\beta -1=0$
and $\beta =2.1149...$. Also, because the discriminant of the polynomial $%
x^{3}-4x-1$ is equal to the prime number $229,$ we see that $\func{disc}(%
\mathbb{Q}(\beta ))=229,$ $U_{\mathbb{Q}(\beta )}=\{1\},$ and so, by Theorem
1.3(vi), the two smallest element of $\tciFourier _{\mathbb{Q}(\beta )}$ are 
$1$ and $\beta .$

\newpage

\textbf{3. Some Lemmas}

In the form in which we use it here, where all numbers are real, Kronecker's
theorem on linearly independent numbers may be stated as follows (see also
[3] and [14, page 66]).

\medskip

\textbf{Lemma 3.0 [8]} \textit{If }$1,$\textit{\ }$\omega _{2},$\textit{\ }$%
...,$\textit{\ }$\omega _{v}$\textit{\ are }$\mathbb{Q}$\textit{-linearly
independent, }$\eta _{2},$\textit{\ }$...,$\textit{\ }$\eta _{v}$\textit{\
are arbitrary, and }$\rho $\textit{\ and }$\varepsilon $\textit{\ are
positive, then there exist integers }$p_{2},$ $...,$\textit{\ }$p_{v},$%
\textit{\ }$n>\rho $\textit{\ such that \ }$\max_{2\leq j\leq v}\left\vert
n\omega _{j}-p_{j}-\eta _{j}\right\vert <\varepsilon .$

\medskip

\bigskip \textbf{Proof.} See for instance [4, Chapter 4].

The first part of the following simple lemma is given in [15, Proposition
3(i)].

\medskip

\textbf{Lemma 3.1 } \textit{Let }$\beta >2$\textit{\ be a real algebraic
integer of degree }$d\geq 2$\textit{\ whose other conjugates lie in the
interval }$(-2,2).$\textit{\ Then, there is a Salem number }$\tau $\textit{\
of degree }$2d$\textit{\ such that }$\beta =\tau +1/\tau ,$\textit{\ i. e., }%
$\beta $ \textit{is a Salem trace number.} \textit{Moreover, }$\mathbb{Q}%
(\beta )$ \textit{is the set of totally real algebraic numbers lying in }$%
\mathbb{Q}(\tau ),$\textit{\ and so it is the unique totally real subfield
of }$\mathbb{Q}(\tau )$\textit{\ with degree }$d.$

\bigskip

\textbf{Proof. }Let $\beta _{1}:=\beta ,$ $\beta _{2},...,\beta _{d}$ be the
conjugates of $\beta ,$ and let $\tau $ be the root greater than $1$ of the
quadratic polynomial $x^{2}-$ $\beta x+1.$ Then, $\tau $ is a zero of the
monic polynomial $P(x):=(x^{2}-\beta _{1}x+1)\cdot \cdot \cdot (x^{2}-\beta
_{d}x+1)\in \mathbb{Z}[x],$ and so it is an algebraic integer of degree at
most $2d,$ whose conjugates are among the two real numbers $\tau $ and $%
1/\tau $ and the $2d-2$ non-real numbers of modulus $1$ roots of the
quadratic polynomials $x^{2}-$ $\beta _{j}x+1,$ where $j\in \{2,...,d\}.$
Further, as $\tau +1/\tau =\beta ,$ we have that $\mathbb{Q}(\beta )\subset 
\mathbb{Q}(\tau ),$ and so $\deg (\tau )\in \{d,$ $2d\}.$ In fact, if $\deg
(\tau )=d,$ then $\mathbb{Q}(\beta )=\mathbb{Q}(\tau ),$ $\tau $ is totally
real with minimal polynomial $x^{2}-$ $\beta x+1,$ and this last assertion
leads immediately to the contradiction $\beta \in \mathbb{Z}.$ Therefore, $%
\deg (\tau )=2d,$ $\tau $ is a Salem number with minimal polynomial $P,$ and 
$\beta $ is a Salem trace number.

To complete the proof of the lemma, notice that any element of $\mathbb{Q}%
(\beta )$ is totally real, and assume on the contrary that there is a
totally real number $\delta \in \mathbb{Q}(\tau )$ which is not in $\mathbb{Q%
}(\beta ).$ Then, the relation $\deg (\mathbb{Q}(\beta ,\delta ))=nd\leq
2d=\deg (K),$ where $n$ is a natural number greater than $1,$ yields $n=2,$ $%
\mathbb{Q}(\tau )$ is equal to the totally real field $\mathbb{Q}(\beta
,\delta ),$ and this last equality contradicts the fact that $\tau $ has
non-real conjugates; thus $\mathbb{Q}(\beta )$ contains all totally real
numbers in $\mathbb{Q}(\tau ).$%
%TCIMACRO{\TeXButton{End Proof}{\endproof}}%
%BeginExpansion
\endproof%
%EndExpansion

\bigskip

The following result is due to Pisot [13] (see also [1] and [14]) for Salem
numbers, and to Boyd (unpublished) for Pisot numbers (see [1, Theorem 8.2],
and [11] for some related generalizations).

\medskip \medskip

\textbf{Lemma 3.2 (Pisot-Boyd)}\textit{\ Let }$e^{\pm i2\pi \omega _{2}},$ $%
...,$ $e^{\pm i2\pi \omega _{d}}$ \textit{be the non-real conjugates of a
Salem number of degree }$2d.$\textit{\ Then, the numbers }$1,$\textit{\ }$%
\omega _{2},$ $...,$ $\omega _{d}$\textit{\ are\ }$\mathbb{Q}$\textit{%
-linearly independent. }

\textit{Similarly, if }$\rho _{1}e^{\pm i2\pi \omega _{1}},$ $...,$ $\rho
_{s}e^{\pm i2\pi \omega _{s}}$ \textit{denote the non-real conjugates of a
non-totally real Pisot number, then the numbers }$1,$\textit{\ }$\omega
_{1}, $ $...,$ $\omega _{s}$\textit{\ are\ }$\mathbb{Q}$\textit{-linearly
independent. }

\medskip

\textbf{Proof.} To give a common proof for the two cases, consider a
non-totally real Pisot number (resp. a Salem number) $\alpha ,$ and let $%
\rho _{2}e^{\pm i2\pi \omega _{2}},$ $...,$ $\rho _{v}e^{\pm i2\pi \omega
_{v}}$ be the non-real conjugates of $\alpha ,$ where $\rho _{2}<1,$ $...,$ $%
\rho _{v}<1$ and $v=s+1\geq 2$ (resp. where $\rho _{2}=\cdot \cdot \cdot
=\rho _{v}=1$ and $v=d).$

Assume on the contrary that there are integers $l,$ $l_{2},...,$ $l_{v}$ not
all zero such that $l_{2}\omega _{2}+\cdot \cdot \cdot +l_{v}\omega _{v}=l.$
By replacing, if necessary, $\omega _{j}$ by $-\omega _{j},$ $\forall $ $%
j\in \{2,$ $...,$ $v\},$ we may assume that all $l_{j}$ are non-negative.
Also, by reordering, if necessary, the non-real conjugates of $\alpha ,$ we
may suppose that $l_{2}=\max \{l_{j}\mid 2\leq j\leq v\}.$ Then, $l_{2}\geq
1,$ $e^{i2\pi l_{2}\omega _{2}}\cdot \cdot \cdot e^{i2\pi l_{v}\omega
_{v}}=e^{i2\pi l}=1,$ $e^{-i2\pi l_{2}\omega _{2}}\cdot \cdot \cdot
e^{-i2\pi l_{v}\omega _{v}}=1,$ $e^{i2\pi l_{2}\omega _{2}}\cdot \cdot \cdot
e^{i2\pi l_{v}\omega _{v}}=e^{-i2\pi l_{2}\omega _{2}}\cdot \cdot \cdot
e^{-i2\pi l_{v}\omega _{v}},$ and 
\[
\alpha _{2}^{l_{2}}\cdot \cdot \cdot \alpha _{v}^{l_{v}}=\alpha
_{v+2}^{l_{2}}\cdot \cdot \cdot \alpha _{2v}^{l_{v}}. 
\]%
where $\alpha _{2}:=\rho _{2}e^{i2\pi \omega _{2}},$ $...,$ $\alpha
_{v}:=\rho _{v}e^{i2\pi \omega _{v}},$ $\alpha _{v+2}:=\rho _{2}e^{-i2\pi
\omega _{2}}=\overline{\alpha _{2}},$ $...,$ $\alpha _{2v}:=\rho
_{v}e^{-i2\pi \omega _{v}}=\overline{\alpha _{v}}.$

Let $\sigma $ be an automorphism of the normal closure $\Gamma $ of $\mathbb{%
Q}(\alpha )$ in $\mathbb{C},$ sending $\alpha _{2}$ to $\alpha .$ Since the
Galois group of $\Gamma $ operates transitively on the conjugates of $\alpha
,$ we have $\left\vert \sigma (\alpha _{v+2})\right\vert <1,$ $...,$ $%
\left\vert \sigma (\alpha _{2v})\right\vert <1$ (resp. $\sigma (\alpha
_{v+2})=\sigma (\overline{\alpha _{2}})=\sigma (1/\alpha _{2})=1/\alpha <1,$ 
$\left\vert \sigma (\alpha _{v+3})\right\vert =1,$ $...,$ $\left\vert \sigma
(\alpha _{2v})\right\vert =1),$ and so 
\[
\alpha ^{l_{2}}\dprod\limits_{j=3}^{v}\left\vert \sigma (\alpha
_{j})\right\vert ^{l_{j}}=\left\vert \sigma (\alpha _{v+2})\right\vert
^{l_{2}}\cdot \cdot \cdot \left\vert \sigma (\alpha _{2v})\right\vert
^{l_{v}}<1. 
\]%
The last relation together with 
\[
1\leq \left\vert \func{Norm}(\alpha )\right\vert ^{l_{2}}\leq \left\vert
\alpha \dprod\limits_{j=3}^{v}\sigma (\alpha _{j})\right\vert ^{l_{2}}\leq
\alpha ^{l_{2}}\dprod\limits_{j=3}^{v}\left\vert \sigma (\alpha
_{j})\right\vert ^{l_{j}}, 
\]%
lead to a contradiction.%
%TCIMACRO{\TeXButton{End Proof}{\endproof}}%
%BeginExpansion
\endproof%
%EndExpansion

\bigskip

The lemma below is also well-known, and may be found in [14] and in [15].
For a seek of completeness we give a proof of this result.

\medskip\ \medskip

\textbf{Lemma 3.3 }\textit{\ Let }$\tau $ \textit{be a Salem number\ of
degree }$2d$ \textit{with non-real conjugates} \textit{\ }$e^{\pm i2\pi
\omega _{2}},$ $...,$ $e^{\pm i2\pi \omega _{d}},$ \textit{and let }$n\in 
\mathbb{N}.$\textit{\ Then, }$\tau ^{n}$\textit{\ is Salem number of degree }%
$2d,$\textit{\ and }$\tau ^{n}+1/\tau ^{n}$ \textit{is a real algebraic
integer greater than }$2$\textit{\ whose other conjugates\ are }$2\cos (2\pi
n\omega _{2}),$ $...,$ $2\cos (2\pi n\omega _{d}),$ \textit{i.e., }$\tau
^{n}+1/\tau ^{n}$\textit{\ is a Salem trace number of degree }$d.$ \textit{%
Moreover, the sequence }$((2\cos (2\pi n\omega _{2}),$ $...,$ $2\cos (2\pi
n\omega _{d}))_{n\in \mathbb{N}}$\textit{\ is dense in} $[-2,2]^{d-1}.$

\medskip

\textbf{Proof. }Let $\sigma _{1},...,\sigma _{2d}$ denote the distinct
embeddings of $\mathbb{Q}(\tau )$ into $\mathbb{C},$ labelled so that $%
\sigma _{1}(\tau )=\tau ,$ $\sigma _{2}(\tau )=e^{i2\pi \omega _{2}},$ $...,$
$\sigma _{d}(\tau )=e^{i2\pi \omega _{d}},$ and $\sigma _{d+j}(\tau
)=1/\sigma _{j}(\tau )$ for all $j\in \{1,$ $...,$ $d\}.$ Then, the
conjugates of $\tau ^{n},$ where $n\in \mathbb{N},$ are the numbers $\sigma
_{1}(\tau )^{n},$ $...,$ $\sigma _{d}(\tau )^{n},$ and as 
\[
\left\vert \sigma _{j}(\tau ^{n})\right\vert \leq 1<\tau ^{n}=\sigma
_{1}(\tau ^{n}),\text{ \ }\forall j\in \{2,\text{ }...,\text{ }2d\}, 
\]%
$\tau ^{n}$ is a Salem number of degree $2d.$ Also, from its definition, a
Salem number is a unit, and hence $\tau ^{n}+1/\tau ^{n}$ is an algebraic
integer. Moreover, as 
\[
\sigma _{1}(\tau ^{n}+1/\tau ^{n})=\sigma _{d+1}(\tau ^{n}+1/\tau ^{n})=\tau
^{n}+1/\tau ^{n}>2 
\]%
and 
\[
-2<\sigma _{j}(\tau ^{n}+1/\tau ^{n})=\sigma _{d+j}(\tau ^{n}+1/\tau
^{n})=2\cos (2\pi n\omega _{j})<2,\text{ \ }\forall j\in \{2,...,d\}, 
\]%
$\tau ^{n}+1/\tau ^{n}$ is repeated twice by the action of the distinct
embeddings of $\mathbb{Q}(\tau )$ into $\mathbb{C};$ thus $\deg (\tau
^{n}+1/\tau ^{n})=d$ and the other conjugates of $\tau ^{n}+1/\tau ^{n}$ are%
\textit{\ }$2\cos (2\pi n\omega _{2}),$ $...,$ $2\cos (2\pi n\omega _{d}).$

To show the last assertion, fix a real number $\varepsilon >0$ and an
element $(r_{2},...,r_{d})$ of $(-2,2)^{d-1}.$ Then, Lemmas 3.2 and 3.0 give
that there are integers $p_{2},$ $...,$ $p_{d},$ $n\geq 1$ such that 
\[
\max_{2\leq j\leq d}\left\vert n\omega _{j}-p_{j}-\frac{\arccos (r_{j}/2)}{%
2\pi }\right\vert <\frac{\varepsilon }{4\pi }. 
\]%
Therefore, for each $j\in \{2,...,d\},$ 
\[
\left\vert \cos (2\pi n\omega _{j}-2\pi p_{j})-\cos (\arccos (\frac{r_{j}}{2}%
))\right\vert <\left\vert (2\pi n\omega _{j}-2\pi p_{j})-\arccos (\frac{r_{j}%
}{2})\right\vert <\frac{\varepsilon }{2},\text{ } 
\]%
$\left\vert 2\cos (2\pi k\omega _{j})-r_{j}\right\vert <\varepsilon ,$ and
so the sequence $((2\cos (2\pi k\omega _{2}),...,2\cos (2\pi k\omega
_{d}))_{k\in \mathbb{N}}$\ is dense in $[-2,2]^{d-1}.$%
%TCIMACRO{\TeXButton{End Proof}{\endproof}}%
%BeginExpansion
\endproof%
%EndExpansion

\bigskip\ 

The proofs of Theorems 1.3(i) and 1.3(iv) are based on the result below.

\medskip

\textbf{Lemma 3.4 } \textit{Let }$K$\textit{\ be a totally real number field
of degree }$d\geq 2.$ \textit{Then, we have the following.}

\textit{(i) If }$\theta \in \wp _{K},$\textit{\ then }$\theta >2$\textit{\
and so }$\theta $ \textit{is a Salem trace number of degree }$d,$ \textit{%
except when} $\theta =(1+\sqrt{5})/2$ (\textit{and} $K=\mathbb{Q}(\sqrt{5}%
)); $

\textit{(ii) the set }$\{(\sigma _{2}(\theta ),...,\sigma _{d}(\theta ))$ $%
\mid \theta \in \wp _{K}\},$\textit{\ where }$\sigma _{2},...,\sigma _{d}$%
\textit{\ are the distinct embeddings of }$K\ $\textit{into }$\mathbb{R}$%
\textit{\ other than the identity of }$K,$\textit{\ is dense in} $%
[-1,1]^{d-1}.$

\medskip\ 

\textbf{Proof.} \textbf{(i)} Let $\theta \in \wp _{K}.$ Then, the other
conjugates of $\theta $ lie in the interval $(-1,1)\subset (-2,2),$ and so $%
\theta $ is a Salem trace number whenever it is greater than $2.$ In fact, a
short calculation shows that $\theta >2^{(d-1)/2},$ as stated by [17, Lemma
2], and so Lemma 3.4(i) is true, since $(1+\sqrt{5})/2$\ is the unique
quadratic Pisot number less than $2.$

\smallskip

\textbf{(ii)} Fix a Salem trace number $\beta ,$ satisfying $\mathbb{Q}%
(\beta )=K.$ Such an element exists by Lemma 3.4(i), since $\wp _{K}$ is
relatively dense in $[1,\infty ).$ Then, the first assertions in Lemmas 3.1
and 3.3 state that there is Salem number $\tau $ of degree $2d$ such that $%
\beta =\tau +1/\tau ,$ and for each $n\in \mathbb{N}$ the number $\beta
_{n}:=\tau ^{n}+1/\tau ^{n}$ is a Salem trace number of degree $d.$ Also,
the last assertions in Lemmas 3.1 and 3.3 give that $\beta _{n}\in \mathbb{Q}%
(\beta )$ and the set $\{(\sigma _{2}(\beta _{n}),...,\sigma _{d}(\beta
_{n}))\mid n\in \mathbb{N}\}$ is dense in $[-2,2]^{d-1}.$ It follows that
the set 
\[
N:=\{n\in \mathbb{N}\mid (\sigma _{2}(\beta _{n}),...,\sigma _{d}(\beta
_{n}))\in (-1,1)^{d-1}\}, 
\]%
is not finite, and the sequence $(\sigma _{2}(\beta _{n}),...,\sigma
_{d}(\beta _{n}))_{n\in N}$ is dense in $[-1,1]^{d-1}.$ Hence, the set $%
\{(\sigma _{2}(\theta ),...,\sigma _{d}(\theta ))$ $\mid \theta \in \wp
_{K}\},$ containing $\{(\sigma _{2}(\beta _{n}),...,\sigma _{d}(\beta
_{n}))\mid n\in N\},$ is also dense in $[-1,1]^{d-1},$ since $\beta _{n}\in
\wp _{K}$ for each $n\in N.$%
%TCIMACRO{\TeXButton{End Proof}{\endproof}}%
%BeginExpansion
\endproof%
%EndExpansion

\bigskip

\textbf{4. Proof of Theorem 1.1 and Corollary 1.2}

As mentioned in the introduction, Theorem 1.1 is true when $K=\mathbb{Q}.$
Suppose that the degree $d$ of the real number field $K$ is greater than $1,$
and let $\sigma _{1},...,\sigma _{d}$ be the distinct embeddings of $K$ into 
$\mathbb{C},$ labelled so that $\sigma _{1}$ is the identity of $K,$ $\sigma
_{1},...,\sigma _{r}$ are real, $\sigma _{r+1},...,\sigma _{r+s}$ are
non-real, and $\sigma _{r+s+1}=c\circ \sigma _{r+1},$ $...,$ $\sigma
_{r+2s}=c\circ \sigma _{r+s},$ where $c$ is the complex conjugation, $\circ $
is the composition function, and $d=r+2s.$

Notice that if $(\alpha ,\alpha ^{\prime })\in \wp _{K}^{2},$ then $\alpha
\alpha ^{\prime }\in \wp _{K},$ since 
\[
\left\vert \sigma _{j}(\alpha \alpha ^{\prime })\right\vert =\left\vert
\sigma _{j}(\alpha )\right\vert \left\vert \sigma _{j}(\alpha ^{\prime
})\right\vert <1<\alpha \alpha ^{\prime }=\sigma _{1}(\alpha \alpha ^{\prime
}),\text{ }\forall \text{ }j\in \{2,...,d\}, 
\]%
and $\alpha \alpha ^{\prime }$ is not repeated under the action of the
distinct embeddings of $K$ into $\mathbb{C};$ thus $\wp _{K}$ is closed
under multiplication, and the powers of any element of $\wp _{K}$ belong to $%
\wp _{K}$ too.

\smallskip

\textbf{(i) }Let $\alpha \in \wp _{K}$ with conjugates $\alpha _{1}:=\sigma
_{1}(\alpha ),...,\alpha _{d}:=\sigma _{d}(\alpha ).$ By replacing, if
necessary, $\alpha $ by $\alpha ^{2}$ we may assume that $\alpha _{1}>0,$ $%
...,$ $\alpha _{r}>0.$ In fact, we shall prove that there are infinitely
many natural numbers $q$ such that

\[
\alpha _{q}^{\prime }:=\alpha ^{q}-1\in \wp _{K}. 
\]%
This implies that $1=\alpha ^{q}-\alpha _{q}^{\prime }\in D_{K},$ and $1$
may be expressed as a difference of two elements of $\wp _{K}$ in infinitely
many ways.

In fact, if $s=0,$ then $\alpha _{q}^{\prime }\in \wp _{K},$ $\forall $ $%
q\in \mathbb{N},$ because the relations $\sigma _{j}(\alpha ^{q}-1)=\alpha
_{j}^{q}-1\in (-1,0),$ \ when $j$ runs through $\{2,...,d\},$ together with
the fact that the norm of $\alpha _{q}^{\prime }$ is a non-zero integer
imply the inequality $\alpha _{q}^{\prime }>1.$

Now suppose $s\geq 1,$ and set $\alpha _{r+1}:=\rho _{1}e^{i2\pi \omega
_{1}},$ $...,$ $\alpha _{r+s}:=\rho _{s}e^{i2\pi \omega _{s}}.$ Then, $%
\omega _{j}\notin \mathbb{Q},$ $\forall $ $j\in \{1,...,s\},$ since
otherwise there are natural numbers $n$ and $j\leq s$ such that $\alpha
_{r+j}^{n}=\rho _{j}^{n}e^{i2\pi n\omega _{j}}=\rho _{j}^{n}\in \mathbb{R},$
and so $\mathbb{Q}(\alpha _{r+j}^{n})\varsubsetneq \mathbb{Q}(\alpha _{r+j})$
contradicting the fact that $\mathbb{Q}(\alpha ^{n})=\mathbb{Q}(\alpha )$
(in fact Lemma 3.2 yields also $\omega _{j}\notin \mathbb{Q},$ $\forall $ $%
j\in \{1,...,s\}).$

It follows by Dirichlet's approximation theorem (see for instance [3]) that
for each real number $Q\geq 6^{s},$ there is a natural number $q\leq Q$ and $%
s$ integers $p_{1},$ $...,$ $p_{s\text{ }}$ such that 
\[
0<\left\vert q\omega _{j}-p_{j}\right\vert <\frac{1}{Q^{1/s}}\leq \frac{1}{6}%
,\text{ }\forall \text{ }j\in \{1,...,s\}. 
\]%
Hence, $\left\vert 2\pi q\omega _{j}-2\pi p_{j}\right\vert <$ $\pi /3,$ $%
\cos (2\pi q\omega _{j})>1/2,$ 
\[
\left\vert \alpha _{r+j}^{q}-1\right\vert ^{2}=(\rho _{j}^{q}e^{i2\pi
q\omega _{j}}-1)(\rho _{j}^{q}e^{-i2\pi q\omega _{j}}-1)=\rho _{j}^{q}(\rho
_{j}^{q}-2\cos (2\pi q\omega _{j}))+1<1, 
\]%
and $\left\vert \alpha _{r+j}^{q}-1\right\vert <1,$ $\forall $ $j\in
\{1,...,s\}.$

Because $\sigma _{j}(\alpha ^{q}-1)=\alpha _{j}^{q}-1$ $\in (-1,0),$ $%
\forall $ $j\in \{2,...,r\},$ whenever $r\geq 2,$ $\left\vert \sigma
_{r+s+j}(\alpha ^{q}-1)\right\vert =\left\vert c(\alpha
_{r+j}^{q}-1)\right\vert =\left\vert \alpha _{r+j}^{q}-1\right\vert <1$ for
all $j\in \{1,...,s\},$ and the norm of the positive algebraic integer $%
\alpha ^{q}-1$ is a non-zero integer, we see that $\alpha ^{q}-1\in \wp _{K}$
and so Theorem 1.1(i) is true. Since $n\omega _{j}-p_{j}\neq 0$ for all $%
(j,n,p_{j})\in \{1,...,s\}\times \mathbb{N}\times \mathbb{Z},$ by letting $Q$
tend to infinity we obtain that there are infinitely many $q\in \mathbb{N}$
such that $\alpha ^{q}-1\in \wp _{K}.$

\bigskip

\textbf{(ii) and Corollary 1.2} \ Let $\tau $ be a Pisot number or a Salem
number, satisfying $K=\mathbb{Q}(\tau ).$ Then, Theorem 1.1(i) gives that $%
1=\alpha -\alpha ^{\prime }$ \ for some $(\alpha ,\alpha ^{\prime })\in \wp
_{K}^{2},$ and so $\tau =\tau \alpha -\tau \alpha ^{\prime }.$ Arguing as
above, we get from \ the relations $\left\vert \sigma _{j}(\tau \alpha
)\right\vert =\left\vert \sigma _{j}(\tau )\sigma _{j}(\alpha )\right\vert
\leq \left\vert \sigma _{j}(\alpha )\right\vert <1$ and $\left\vert \sigma
_{j}(\tau \alpha ^{\prime })\right\vert =\left\vert \sigma _{j}(\tau )\sigma
_{j}(\alpha ^{\prime })\right\vert \leq \left\vert \sigma _{j}(\alpha
^{\prime })\right\vert <1,$ $\forall $ $j$ $\in \{2,...,d\},$ that the
positive algebraic integers $\tau \alpha $ and $\tau \alpha ^{\prime }$
belong to $\wp _{K},$ and so $\tau \in D_{K}.$ In particular, $\wp
_{K}\subset D_{K},$ and $\wp _{K}\varsubsetneq D_{K},$ as $1\in D_{K}.$

\smallskip

\textbf{(iii)} The proof of the second assertion is immediate. Indeed, if $%
\func{card}(\tciFourier _{K})=1,$ then (2) gives that for any $n\in \mathbb{N%
},$ $\theta _{n}-\theta _{1}=(n-1)d_{1},$ 
\[
(n-1)\left\vert \sigma _{j}(d_{1})\right\vert =\left\vert \sigma _{j}(\theta
_{n}-\theta _{1})\right\vert \leq \left\vert \sigma _{j}(\theta
_{n})\right\vert +\left\vert \sigma _{j}(\theta _{1})\right\vert \leq 2, 
\]%
where $j$ is any fixed element of $\{2,...,d\},$ and this last relation is
not true when $n$ sufficiently is large.

To show the first assertion, recall from the proof of Theorem 1.1(i) that if 
$\alpha \in \wp _{K}$ and the real conjugates of $\alpha $ are positive,
then there are infinitely many natural numbers $q$ such that $\alpha
^{q}-1\in \wp _{K}.$ If we apply this result, when $r\geq 2,$ to the square
of an element $\alpha $ of $\wp _{K},$ we see that each element of $D_{K}$
of the form $(\alpha ^{2q}-1)-\alpha ^{q}$ does belong to $\wp _{K},$ as $%
\sigma _{2}(\alpha )^{2q}-1-\sigma _{2}(\alpha )^{q}<-1,$ and so Theorem
1.1(iii) is true, whenever $r\geq 2.$

Now, suppose that all other conjugates of $\alpha $ are non-real, and we
claim there are infinitely many natural numbers $n$ such that $\alpha
^{n}+1\in \wp _{K}.$ Indeed, setting $\sigma _{2}(\alpha ):=\rho
_{1}e^{i2\pi \omega _{1}},$ $...,$ $\sigma _{1+s}(\alpha ):=\rho
_{s}e^{i2\pi \omega _{s}},$ we have by Lemma 3.2 that $1,$\textit{\ }$%
2\omega _{1},$ $...,$ $2\omega _{s}$\textit{\ }are\textit{\ }$\mathbb{Q}$%
-linearly independent, and so it follows by Lemma 3.0 that for each $%
\varepsilon \in (0,1/6)$ there are integers $p_{1},$ $...,$ $p_{s}$ and a
natural number $n$ such that \ 
\[
\max_{1\leq j\leq s}\left\vert n\omega _{j}-p_{j}-\frac{1}{2}\right\vert
<\varepsilon <\frac{1}{6}. 
\]%
Hence, for each $j\in \{1,...,s\},$ $\left\vert 2n\pi \omega _{j}-2\pi
p_{j}-\pi \right\vert <\pi /3,$ $\cos (2\pi n\omega _{j})<-1/2,$ and so 
\[
\left\vert \sigma _{1+j}(\alpha )^{n}+1\right\vert ^{2}=(\rho
_{j}^{n}e^{i2\pi n\omega _{j}}+1)(\rho _{j}^{n}e^{-i2\pi n\omega
_{j}}+1)=\rho _{j}^{q}(\rho _{j}^{q}+2\cos (2\pi q\omega _{j}))+1<1, 
\]%
and in the same as in the end of the proof of Theorem 1.1(i), we get that $%
\alpha ^{n}+1\in \wp _{K}.$ Also, since $n\omega _{j}-p_{j}-\frac{1}{2}\neq
0 $ for all $(j,n,p_{j})\in \{1,...,s\}\times \mathbb{N}\times \mathbb{Z},$
by letting $\varepsilon $ tend to zero we get that there are infinitely many 
$n$ such that $\alpha ^{n}+1\in \wp _{K}$ and this ends the proof of the
claim.

We can now easily deduce that Theorem 1.1(iii) is true. Indeed, choose a
large integer $q$ so that $\alpha ^{q}-1\in \wp _{K}$ and $\left\vert \sigma
_{2}(\alpha )\right\vert ^{q}<1/2,$ and let $n$ be any integer satisfying $%
\alpha ^{n}+1\in \wp _{K}$ and $n>q.$ Then, any element $\beta _{n}:=(\alpha
^{n}+1)-(\alpha ^{q}-1)$ of $D_{K},$ satisfies $\left\vert \sigma _{2}(\beta
_{n})\right\vert =\left\vert 2+\alpha ^{n}-\alpha ^{q}\right\vert
>2-\left\vert \sigma _{2}(\alpha )\right\vert ^{n}-\left\vert \sigma
_{2}(\alpha )\right\vert ^{q}>2-2\left\vert \sigma _{2}(\alpha )\right\vert
^{q}>1,$ and so $\beta _{n}\notin \wp _{K}.$%
%TCIMACRO{\TeXButton{End Proof}{\endproof}}%
%BeginExpansion
\endproof%
%EndExpansion

\bigskip \bigskip

\textbf{5. Proof of Theorem 1.3}

Let $K$ be totally real number field of degree $d\geq 2,$ and let $\sigma
_{1},$ $...,$ $\sigma _{d}$ be the distinct embeddings of $K$ into $\mathbb{R%
},$ where $\sigma _{1}$ is the identity of $K.$

\textbf{(i)} As mentioned in the relation (1), $D_{K}$ is always contained
in $E_{K}.$ To show the inclusion $E_{K}\subset $ $D_{K},$ consider an
element $\beta $ of $E_{K}.$ Then, $\beta $ is a positive algebraic integer
of $K$ such that $(\sigma _{2}(\beta ),...,\sigma _{d}(\beta ))\in
(-2,2)^{d-1}.$ Recall, by Lemma 3.4(ii), that the set $\{(\sigma _{2}(\theta
),...,\sigma _{d}(\theta ))$ $\mid \theta \in \wp _{K}\}$\textit{\ }is dense
in $[-1,1]^{d-1},$ and hence we can find (infinitely many) $\theta \in \wp
_{K}$ such that for each $j\in \{2,...,d\},$ 
\[
-1<\sigma _{j}(\theta )<1-\sigma _{j}(\beta )\text{ \ if }\sigma _{j}(\beta
)>0, 
\]%
and 
\[
-1-\sigma _{j}(\beta )<\sigma _{j}(\theta )<1\text{\ \ when }\sigma
_{j}(\beta )<0. 
\]%
It follows that the real algebraic integer $\theta ^{\prime }:=\beta +\theta
>1$ belongs to $K$ and satisfies $\sigma _{j}(\theta ^{\prime })=\sigma
_{j}(\beta )+\sigma _{j}(\theta )\in (-1,1)$ for all $j\in \{2,...,d\};$
thus $\deg (\theta ^{\prime })=d,$ as it is not repeated under the action of
the distinct embeddings of $K$ into $\mathbb{R},$ $\theta ^{\prime }\in \wp
_{K}$ and $\beta =\theta ^{\prime }-\theta \in D_{K}.$

\smallskip

\textbf{(ii)} From the definition of a root of unity trace number we have
that $U_{K}\subset E_{K}\cap (0,2).$ To prove the inclusion inverse consider
an algebraic integer $\beta $ of $K$ such that $(\sigma _{1}(\beta ),\sigma
_{2}(\beta ),...,\sigma _{d}(\beta ))\in (0,2)\times (-2,2)\times \cdot
\cdot \cdot \times (-2,2).$ Because the zeros of the monic polynomial 
\[
P(x):=(x^{2}-\sigma _{1}(\beta )x+1)\cdot \cdot \cdot (x^{2}-\sigma
_{d}(\beta )x+1)\in \mathbb{Z}[x] 
\]
are all non-real numbers of modulus $1,$ we have, by Kronecker's theorem
(see for instance [3]), that these zeros are roots of unity. Hence, there
are natural numbers $n$ and $k\leq n$ with $\gcd (k,n)=1$ such that $\beta
=\sigma _{1}(\beta )=e^{i2\pi k/n}+1/e^{i2\pi k/n}=2\cos (2k\pi /n).$
Further, the condition $\beta \in (0,2)$ implies that $n\geq 5,$ $k\leq n-1$
and $\beta \in U_{K};$ thus $U_{K}=E_{K}\cap (0,2),$ and so, by Theorem
3.1(i), $U_{K}=D_{K}\cap (0,2).$

\smallskip

\textbf{(iii) } The first assertion follows from Theorems 1.3(i) and
1.3(ii). It is also clear that any element of the set $\{\beta \in \mathbb{Z}%
_{K}\cap (2,\infty )\mid (\sigma _{2}(\beta ),...,\sigma _{d}(\beta ))\in
(-2,2)^{d-1}\}$ is not repeated under the action of the distinct embeddings
of $K$ into $\mathbb{R}$ and so it is of degree $d.$

\smallskip

\textbf{(iv)} Fix an element $(\varepsilon _{2},...,\varepsilon _{d})$ \ of
the set $\{-1,1\}^{d-1}.$ Then, Lemma 3.4(ii) gives that there is a
subsequence $(\theta _{f(n)})_{n\in \mathbb{N}}$ of the sequence $(\theta
_{n})_{n\in \mathbb{N}}$ such that $\lim_{n\rightarrow \infty }(\sigma
_{2}(\theta _{f(n)}),...,\sigma _{d}(\theta _{f(n)}))=(\varepsilon
_{2},...,\varepsilon _{d}).$ Because $\{\theta _{f(n)+1}-\theta _{f(n)}\mid $
$n\in \mathbb{N\}}$ is contained in the finite set $\tciFourier _{K},$ there
exist an element $d_{\ast }$ of $\tciFourier _{K}$ and an infinite subset $N$
of $\mathbb{N}$ such that $\theta _{f(n)+1}-\theta _{f(n)}=d_{\ast }$ \ for
all $n\in N.$ Therefore for each $(j,n)\in \{2,...,d\}\times N,$ $\sigma
_{j}(\theta _{f(n)})+\sigma _{j}(d_{\ast })=\sigma _{j}(\theta _{f(n)+1}),$
and hence 
\[
-1<\sigma _{j}(\theta _{f(n)})+\sigma _{j}(d_{\ast })<1. 
\]%
Letting the element $n$ of $N$ tend to infinity, we get from the last
relation that $-1-\varepsilon _{j}\leq \sigma _{j}(d_{\ast })\leq
1-\varepsilon _{j},$ and so $\sigma _{j}(d_{\ast })<0$ (resp. $\sigma
_{j}(d_{\ast })>0)$ when $\varepsilon _{j}=1$ (resp. $\varepsilon _{j}=-1).$

Therefore, for any $(\varepsilon _{2},...,\varepsilon _{d})\in
\{-1,1\}^{d-1},$ there is $d_{\ast }\in \tciFourier _{K}$ such that each $%
(j,n)\in \{2,...,d\},$ the numbers $\sigma _{j}(\beta )$ and $\varepsilon
_{j}$ have opposite signs, and so there are $2^{d-1}$ distinct elements of
the set $\tciFourier _{K},$ which are in a one-to-one correspondence with
the elements $\{-1,1\}^{d-1}.$

\smallskip

\textbf{(v)} Since $\tciFourier _{K}\cap (0,2)\subset D_{K}\cap (0,2),$ the
first assertion follows from Theorem 1.3(ii). Also, Theorem 1.3(ii) says
that small elements of $D_{K}$ belong to $U_{K}.$ In particular, $\min
D_{K}=\min U_{K}=$ $u_{1},$ and so $u_{1}\leq \min \{1,d_{1}\},$ since $%
\{d_{1},1\}\subset D_{K}.$ Using (2), we get $u_{1}=l_{1}d_{1}+\cdot \cdot
\cdot +l_{\func{card}(\tciFourier _{K})}d_{\func{card}(\tciFourier
_{K})}\geq (l_{1}+\cdot \cdot \cdot +l_{\func{card}(\tciFourier
_{K})})d_{1}\geq d_{1},$ for some non-negative integers $l_{1},$ $...,$ $l_{%
\func{card}(\tciFourier _{K})},$ and hence $u_{1}=d_{1}\leq 1.$

Now, suppose $Card(U_{K})\geq 2.$ Then, $\min (D_{K}\backslash \{u_{1}\})=$ $%
u_{2},$ and $u_{2}\leq d_{2},$ as $d_{2}\in D_{K}\backslash
\{d_{1}\}=D_{K}\backslash \{u_{1}\}.$ Also, writing $u_{2}=l_{1}d_{1}+\cdot
\cdot \cdot +l_{\func{card}(\tciFourier _{K})}d_{\func{card}(\tciFourier
_{K})},$ for some non-negative integers $l_{1},$ $...,$ $l_{\func{card}%
(\tciFourier _{K})},$ we obtain by the relation $u_{2}\leq d_{2}$ that $%
u_{2}=d_{2}$ or $u_{2}=l_{1}d_{1}$ from some natural $l_{1}\geq 2.$ In fact
this last possibility does not hold because the equalities $%
u_{2}=l_{1}d_{1}=l_{1}u_{1}$ imply that the conjugates of the non-zero
algebraic integer $u_{1}=u_{2}/l_{1}$ belong to $(-1,1);$ thus $u_{2}=d_{2}$
(in the same way we can show that $u_{3}=d_{3},$ whenever $Card(U_{K})\geq 3$
and $u_{3}$ is not a linear combination of $u_{1}$ and $u_{2}$ with
non-negative integer coefficients).

\smallskip

\textbf{(vi)} Suppose $Card(U_{K})=1.$ Then,\textit{\ }$U_{K}=\{1\},$\textit{%
\ }as $1\in U_{K}.$ It follows by Theorems 1.3(v) and 1.3(iii),
respectively, that \textit{\ }$d_{1}=1$\textit{\ }and $d_{2}\in T_{K};$ thus 
$\min T_{K}\leq d_{2}.$ In the same way as in the last proof, we get $\min
T_{K}=d_{2}$ or $\min T_{K}=l_{1}d_{1}=l_{1}\in \mathbb{N}.$ Because the
degree of each element of $T_{K}$ is equal to $d,$ as asserted by Theorem
1.3(iii), this last possibility does not hold, and so $\min T_{K}=d_{2}.$%
%TCIMACRO{\TeXButton{End Proof}{\endproof}}%
%BeginExpansion
\endproof%
%EndExpansion

\bigskip \bigskip

\textit{In a private communication and after having finished writing this
manuscript, I was informed by Dubickas that he considered in }[6]\textit{\
the more general problem about algebraic integers which may written as a
difference of two Pisot numbers. Using another approach, he obtained
certains results improving some ones of the present manuscript. }

\textit{\bigskip }

\textbf{References}

[1] M. J. Bertin, A. Decomps-Guilloux, M. Grandet-Hugo, M.
Pathiaux-Delefosse and J. P. Schreiber, \textit{Pisot and Salem numbers},
Birkh\"{a}user Verlag Basel, 1992.

[2] M. J. Bertin and T. Za\"{\i}mi, \textit{Complexes Pisot numbers in
algebraic number fields}, C. R. Math., Acad. Sci. Paris \textbf{353}, No. 11
(2015), 965-967.

[3] Y. Bugeaud, \textit{Distribution modulo one and diophantine approximation%
}, Cambridge university press, 2012.

[4] J. W. S. Cassels, \textit{An introduction to Diophantine Approximation,}
Cambridge University Press, 1952.

[5] A. Dubickas, \textit{Every Salem number is a difference of two Pisot
numbers}, Proc. Edinb. Math. Soc., II. Ser. \textbf{66}, No. 3 (2023),
862-867.

[6] A. Dubickas, \textit{Every Salem number is a difference of two Pisot
numbers}, (preprint).

[7] A. H. Fan and J. Schmeling, $\varepsilon -$\textit{Pisot numbers in any
real algebraic number field are relatively dense}, J. Algebra \textbf{272}
(2004), 470-475.

[8] J. C. Lagarias, \textit{Meyer's concept of quasicrystal and quasiregular
sets}, Commun. Math. Phys. \textbf{179} (1996), 365-376.

[9] L. Kronecker, Zwei S\"{a}tze \"{u}ber Gleichungen mit ganzzahligen
Coefficienten, J. Reine Angew. Math., \textbf{53} (1857), 173-175

[10] Y. Meyer, \textit{Nombres de Pisot, nombres de Salem et analyse
harmonique,} Lecture Notes in Mathematics, Springer-Verlag, \textbf{117}
(1970).

[11] M. Mignotte, \textit{sur les conjugu\'{e}s des nombres de Pisot, } C.
R. Acad. Sci. Paris, \textbf{298} (1984), 21-24.

[12] R. V. Moody, \textit{Meyer sets and their duals,} The Mathematics of
long-range aperiodic order (Waterloo, ON, 1995), 403-441, NATO Adv. Sci.
Inst. Ser. C Math. Phys. Sci., \textbf{489}, Kluwer Acad. Publ., Dordrecht,
1997.

[13] C. Pisot, \textit{Quelques aspects de la th\'{e}orie des entiers alg%
\'{e}briques,} S\'{e}minaire de math\'{e}matiques sup\'{e}rieures, Universit%
\'{e} de Montr\'{e}al, 1963.

[14] R.\ Salem, \textit{Algebraic numbers and Fourier analysis,} D. C. Heath
and Co., Boston, Mass., 1963.

[15] C. J. Smyth, \textit{Seventy years of Salem numbers,} Bull. Lond. Math.
Soc. 47, No. 3, 379-395 (2015).

[16] THE PARI GROUP, \textit{PARI/GP version} 2.5.1, 2012, avaible from

http://pari.math.u-bordeaux.fr/.

[17] T. Za\"{\i}mi, \textit{Sur les nombres de Pisot totalement reels,} Arab
J. Math. Sci \textbf{5, }No. 2 (1999)\textbf{, }19-32.

\smallskip \lbrack 18] T. Za\"{\i}mi, \textit{On }$\varepsilon -$\textit{%
Pisot numbers}, New York J. Math. \textbf{15 }(2009)\textbf{, }415-422.

[19] T. Za\"{\i}mi, \textit{Commentaires sur quelques r\'{e}sultats sur les
nombres de Pisot}, J. Th\'{e}or. Nombres Bordx. \textbf{22}, No. 2 (2010),
513-524.

[20] T. Za\"{\i}mi, \textit{Comments on Salem polynomials,} Arch. Math. 
\textbf{117} (2021), 41--51.

\bigskip

\bigskip

Department of Mathematics and Statistics \ College of Science

Imam Mohammad Ibn Saud Islamic University (IMSIU)

P. O. Box 90950 \ 

Riyadh 11623 Saudi Arabia

Email: tmzaemi@imamu.edu.sa\textit{\ }

\bigskip\ 

\bigskip\ 

\bigskip

\textit{\ }

\end{document}